\documentclass[freqno,11pt]{amsart}
\usepackage{amscd,amsfonts,amsmath,amssymb,amstext,amsthm}

\usepackage[utf8]{inputenc}
\usepackage[T2A]{fontenc}
\usepackage[russian,english,]{babel}
\usepackage{cite}
\usepackage[unicode]{hyperref}
\usepackage{color}
\usepackage{xcolor}
\usepackage{comment}
\usepackage[normalem]{ulem}

\makeatletter
\def\@settitle{\begin{center}%
    \baselineskip14\p@\relax
    \bfseries
    \MakeUppercase{\@title}
  \end{center}
}
\makeatother

\newtheorem{theorem}{Theorem}
\newtheorem{lemma}{\it Lemma}[section]
\newtheorem{proposition}{Proposition}[section]
\newtheorem{corollary}{\it Corollary}[section]
\theoremstyle{remark}

\newtheorem{remark}{Remark}[section]
\newtheorem{example}{\it Example}[section]
\theoremstyle{definition}
\newtheorem{assert}{Assertion}
\newtheorem{definition}{\it Definiton}[section]
\numberwithin{equation}{section}

\def\C{{\mathbb C}}
\def\Q{{\mathbb Q}}
\def\Cn*{{\C^n}^*}
\def\CN*{{\C^N}^*}

\def\P{{\mathbb P}}

\def\R{{\mathbb R}}
\def\Z{{\mathbb Z}}

\def\tr{{\rm Tr}}

\def\vol{{\rm vol}}
\def\E{\:{\rm e}}
\def\re{{\rm Re}}

\begin{document}
\title{On real roots of polynomials in the context of group theory
}
{
\author{B. Kazarnovskii
}
\address {\noindent
Moscow Institute of Physics and
Technology (National Research University),
Higher School of Contemporary Mathematics
\newline
{\it kazbori@gmail.com}.}
}
\thanks {\noindent
This research was carried out at the Higher School of Contemporary Mathematics of Moscow
Institute of Physics and Technology,
with the support of the Ministry of Science and Higher Education
of the Russian Federation,
project no. SMG - 2024 -0048.
}

\begin{abstract}
The probability that a root of a random real polynomial of increasing degree is real tends to zero.
The transition from polynomials to Laurent polynomials leads to an unexpected result:
the probability of the root being real tends not to zero but to $1/\sqrt 3$.
A similar phenomenon was also described for systems of $n$ Laurent polynomials in $n$ variables.
Considering Laurent polynomials as functions associated with representations of a torus,
we describe a similar phenomenon for representations of any reductive linear group.
If the group is simple,
a formula for the mentioned probability limit is provided.
\end{abstract}
\maketitle
%

\section{Introduction}\label{s1}
Let the coefficients of a random real polynomial of degree $m$
in one variable be normally and independently distributed
with zero mean and unit variance.
Denote by $\mathcal P(m)$ the probability that a root of the polynomial is real.
Then, if $m$ is large, we have $\mathcal P(m)\asymp\frac{2}{\pi}\frac{\log m}{m}$; see \cite{KA}.
For more details on the distribution of the number of real roots of polynomials,
see the survey \cite{EK} and the bibliography therein.

\subsection{Laurent polynomials}\label{s1.1}
The transition from polynomials to Laurent polynomials leads to an unexpected result:
the probability that a root of a Laurent polynomial of increasing degree is real tends not to zero but to $1/\sqrt 3$;
see Corollary \ref{corAntiKac} or Example \ref{exKac}.
This result is equivalent to computing the asymptotics of the mean number
of zeros of trigonometric polynomials of increasing degree on the circle; see \cite{ADG}.
The phenomenon of a non-zero limiting probability persists
for Laurent polynomials in several variables.
We present here the formula for computing this probability from \cite{K22}.

Recall that
a Laurent polynomial is a function on the complex torus $(\C\setminus0)^n$ of the form
$$P(z)=\sum_{m\in\Lambda\subset\Z^n,\: a_m\in\C}a_m z^m,$$
where $z^m=z_1^{m_1}\cdots z_n^{m_n}$.
A finite set $\Lambda$ in $\Z^n$ is called the \emph{support of the polynomial} $P$.
\begin{definition}\label{dfRealLaur}
A polynomial $P$ is called a \emph{real Laurent polynomial}
if its values on the subtorus
$
T^n=\{z\in(\C\setminus0)^n\colon\:z=(\E^{i\theta_1},\ldots,\E^{i\theta_n})\}
$
are real.
Any root of a Laurent polynomial contained in $T^n$
is called a \emph{real root of the polynomial}.
\end{definition}
\begin{corollary}\label{corReal}
{\rm(1)}\ A Laurent polynomial $\sum_ka_kz^k$ is real
if and only if $\forall k\in\Z^n\colon a_k=\overline{a_{-k}}$.
In particular,
the support of a real Laurent polynomial is centrally symmetric.

{\rm(2)}\
The set of roots of a real Laurent polynomial is invariant under the mapping
$z\mapsto\bar z^{-1}$.
\end{corollary}
Denote by $\mathcal P(\Lambda)$ the probability that a root of a system of $n$ random real Laurent polynomials with support $\Lambda$ is real
(the precise definition of randomness will be given below in \textsection\ref{s2.1}).

Let $B_m$ be the ball in $\R^n$ of radius $m$ centered at the origin,
and $\Lambda_m=B_m\cap \Z^n$, where $\Z^n$ is the integer lattice in $\R^n$.
Laurent polynomials with supports in $\Lambda_m$ can be viewed as multidimensional analogs of Laurent polynomials of degree $m$.
In \cite{K22} it is proved that
\begin{equation}\label{mainLaur}
\lim_{m\to\infty} \mathcal P(\Lambda_m)=\left(\frac{\sigma_{n-1}}{\sigma_n}\beta_n\right)^{\frac{n}{2}}
\end{equation}
where $\beta_n=\int_{-1}^1x^2(1-x^2)^{\frac{n-1}{2}}dx$, and $\sigma_k$ is the volume of the $k$-dimensional ball of radius $1$.
We list
the values of $\beta_n$ for $1\leq n\leq 10$:

\begin{center}
$\beta_n=\frac{2}{3},\:\frac{\pi}{8},\:\frac{4}{15},\:
     \frac{\pi}{16},\:
\frac{16}{105},\:
\frac{5\pi}{128},\:
\frac{32}{315},\:
\frac{7\pi}{256},\:
\frac{256}{3465},\:
\frac{21\pi}{1024}$.
\end{center}
\begin{corollary}\label{corAntiKac}
If $n=1$, then $\lim_{m\to\infty} \mathcal P(\Lambda_m)=1/\sqrt 3$.
\end{corollary}
\begin{proof}
Since $\beta_1=2/3$, we have
$\lim_{m\to\infty} \mathcal P(\Lambda_m)=\left(\frac{\sigma_0}{\sigma_1} \frac{2}{3}\right)^{\frac{1}{2}}=1/\sqrt 3$.
\end{proof}
\subsection{Description of the results}\label{s1.2}
Let $\pi$ be a finite-dimensional representation of a Lie group $G$,
and let ${\rm Trig}(\pi)$ denote the vector space of functions on $G$
consisting of linear combinations of matrix elements of the representation.
If the representation $\pi$ is real (resp. complex),
then ${\rm Trig}(\pi)$ is regarded as a real (resp. complex) vector space.
Functions from the space ${\rm Trig}(\pi)$ are called \emph{$\pi$-polynomials} on $G$.

Let $K$ be a compact connected group.
Denote by $K^\C$ the complexification of the group $K$.
Recall that the group $K^\C$ exists, is unique, and is determined by the following conditions:

 (i)\ $K^\C$ is a connected complex Lie group, $\dim_\C(K^\C)=\dim(K)$;

 (ii)\ the Lie algebra of $K^\C$ is the complexification of the Lie algebra of $K$;

 (iii)\ $K$ is a maximal compact subgroup of $K^\C$.

For example, $(\C\setminus 0)^n$ and $GL(n,\C)$ are complexifications
of the torus $T^n$ and the group of unitary matrices $U(n,\C)$, respectively.

Any complex (resp. real) representation $\pi$ of the group $K$ in a complex (resp. real) space $E$
extends uniquely to its complexification,
i.e., to a holomorphic representation of the group $K^\C$ in $E$ (resp. in $E\otimes_\R\C$).
Hence, any $\pi$-polynomial can also be regarded as a $\pi^\C$-polynomial on the group $K^\C$.
A root lying in $K$ is called a \emph{real root of a $\pi$-polynomial}.

In the context of Laurent polynomials, the notion of a $\pi$-polynomial means the following.
Let $\pi_0$ be the trivial representation of the torus $T^n$ in $\R^1$.
For $0\ne m\in\Z^n$ set
\[
\pi_m(\theta)=\begin{pmatrix}
 \cos(m,\theta) & \sin(m,\theta)\\
 -\sin(m,\theta) & \cos(m,\theta)
\end{pmatrix}.
\]
For any unordered pair $(m,-m)$,
set $\pi_{(m,-m)}=\pi_m$.
Since the real representations $\pi_m$ and $\pi_{-m}$ of the group $T^n$ are equivalent,
the notation $\pi_{(m,-m)}$ is well-defined.
For any finite centrally symmetric set $\Lambda\subset\Z^n$,
denote by $\Lambda'$ the set of unordered pairs $(m,-m)$ with $m\in\Lambda$.
Consider the real representation of the torus $T^n$
\begin{equation}\label{eqpi1}
  \pi(\Lambda)=\bigoplus_{(m,-m)\in\Lambda'}\pi_{(m,-m)}.
\end{equation}
Then

(1) $\pi(\Lambda)$-polynomials on $T^n$ are trigonometric polynomials
\[
f(\E^{i\theta})=\sum_{m\in \Lambda,\:  \alpha_m,\beta_m\in\R} \alpha_m\cos(m,\theta)+\beta_m\sin(m,\theta),
\]

(2) the space of Laurent polynomials with support in $\Lambda$ is
the space of $(\pi(\Lambda))^\C$-polynomials ${\rm Trig}((\pi(\Lambda))^\C)$;

(3) holomorphic extensions of $\pi(\Lambda)$-polynomials to $(\C\setminus0)^n$
are precisely real Laurent polynomials with support $\Lambda$;
see Corollary~\ref{corReal}.

In terms of $\pi$-polynomials on the torus $T^n$,
formula~(\ref{mainLaur}) takes the following form.
Let $\Lambda_m=\Z^n\cap B_m$,
where $B_m$ is the ball of radius $m$ centered at the origin.
Set $\pi(m)=\oplus_{k\in\Lambda_m}\pi_{k,-k}$.
Then the space of $\pi(m)$-polynomials,
extended to Laurent polynomials on $(\C\setminus0)^n$,
coincides with the space of real Laurent polynomials of degree $\leq m$,
and $\mathcal P(\Lambda_m)$ is the probability
that a root of a system of $\pi(m)$-polynomials is real.
The proof of (\ref{mainLaur}) relies on the notions
of the Newton ellipsoid and the Newton polytope of a representation of a torus.
Recall that the Newton polytope of a Laurent polynomial
is the convex hull ${\rm conv}(\Lambda)$
of its support $\Lambda$.
We also associate with the support $\Lambda$ an ellipsoid ${\rm Ell}(\Lambda)$,
called its \emph{Newton ellipsoid} (see Definition~\ref{dfSuppRepr}).
If the support $\Lambda$ is centrally symmetric,
then ${\rm Ell}(\Lambda)\subset{\rm conv}(\Lambda)$.

In \cite{K22} it was proved that the average number of real roots of a system of $n$ real Laurent polynomials with support $\Lambda$
is equal to the volume of the Newton ellipsoid ${\rm Ell}(\Lambda)$.
From this, using Kushnirenko's theorem on the number of roots of a polynomial system,
we obtain
\begin{equation}\label{eqs1.2-1}
 \mathcal P(\pi(\Lambda))=\frac{\vol\left({\rm Ell}(\Lambda)\right)}{\vol\left({\rm conv}(\Lambda)\right)}.
\end{equation}
The proof of (\ref{mainLaur}) is based on (\ref{eqs1.2-1}).
For Laurent polynomials with different supports, a similar equality of the form (\ref{eqs1.2-1}) is valid,
where the volumes in numerator and denominator
are replaced by mixed volumes of the corresponding Newton ellipsoids and Newton polytopes.
In connection with the problem of roots of a random system of equations,
mixed volumes of ellipsoids first appeared in \cite{ZK}.

Below we replace the torus $T^n$ and its real representation
by an arbitrary $n$-dimensional compact group $K$ and its real representation $\pi$.
In this setting, we define the notion of the \emph{Newton ellipsoid of a representation $\pi$};
see Definition~\ref{dfSuppRepr}.
The Newton ellipsoid is an ellipsoid invariant under the coadjoint action of the group,
lying in the space of linear functionals on the Lie algebra of $K$.
We prove that the average number of real roots of a system of $\pi$-polynomials
equals the volume of the Newton ellipsoid multiplied by $\frac{n!}{(2\pi)^n}$;
see Theorem~\ref{thmMixed}.
The proof of this statement is based on the formula for the average number of common zeros of $n$ smooth functions on
an $n$-dimensional differentiable manifold from \cite{AK}; see also \cite{K20}.
\begin{example}\label{exKac}
Consider the representation $\pi=\oplus_{0\leq k\leq m} \pi_{(k,-k)}$ of the torus $T^1$.
In this case, the Newton ellipsoid ${\rm Ell}(\pi)$
is the segment $[-\alpha,\alpha]$,
where, according to Definition~\ref{dfSuppRepr},
$\alpha=2\pi\sqrt{\frac{2}{2m+1}\sum_{-m\leq k\leq m} k^2}=2\pi\sqrt{\frac{m(m+1)}{3}}$.
Hence it follows
that the average number of zeros of trigonometric polynomials in one variable of degree $m$ is $2\sqrt{\frac{m(m+1)}{3}}$.
These trigonometric polynomials are restrictions to $T^1$ of real Laurent polynomials of degree $m$.
Thus, the probability that a root of a real Laurent polynomial of degree $m$ is real equals
$\sqrt{\frac{m+1}{3m}}$.
In particular, as $m\to\infty$
this probability tends to $\sqrt{1/3}$.
\end{example}

In \textsection\ref{BKK}
we define the notion of the \emph{Newton body $\mathfrak N(\mu)$ of a complex representation} $\mu$ of the group $K^\C$;
see Definition~\ref{dfweightedpolyhedron}.
The Newton body is an analogue of the Newton polytope of a Laurent polynomial.
It is a compact convex set in the space of linear functionals on the Lie algebra of $K$,
invariant under the coadjoint action of the group.

Using the above notions of the Newton ellipsoid and the Newton body of a representation,
we prove an analogue of formula (\ref{eqs1.2-1}) for the probability $\mathcal P(\pi)$ of a root of a system of $\pi$-polynomials being real
(see Theorem~\ref{thmProportion}):
\begin{equation}\label{eqPr}
 \mathcal P(\pi)=c(K)\frac{\vol({\rm Ell}(\pi))}{\vol(\mathfrak N(\pi^\C))},
\end{equation}
where $c(K)$ is a constant depending only on the group $K$.
Equality (\ref{eqPr}) remains valid for the case of $\pi_i$-polynomials of different representations $\pi_i$,
if the volumes in numerator and denominator
are replaced by mixed volumes of the corresponding Newton ellipsoids and Newton bodies.

The proof of (\ref{eqPr}) uses
a formula of Kushnirenko–Bernstein–Khovanskii type (also called the BKK formula) for complex reductive groups; see, e.g., \cite{K87,Br,AKK,KK1}.
The version of this formula obtained here (Theorem~\ref{thmRed} in \textsection\ref{BKK})
is closer to the standard formulation of the BKK formula than previously known versions.
More precisely,
we prove that for any representations $\mu_1,\ldots,\mu_n$,
the number of common zeros of a generic system of $n$ functions $f_i\in{\rm Trig}(\mu_i^\C)$,
up to a constant factor depending on the group $K$,
is equal to the mixed volume of the Newton bodies $\mathfrak N(\mu_1),\ldots,\mathfrak N(\mu_n)$.
The proof of this statement is given in \textsection\ref{proportion}; see Theorem~\ref{thmRed1}.

If the group $K$ is simple,
then, using (\ref{eqPr}),
we obtain a formula for the limit probability $\mathcal P(\pi_m)$ for a growing sequence of representations $\pi_m$,
analogous to formula (\ref{mainLaur});
see Theorem~\ref{thmAntiKac2}.
\par\medskip
In Section~\ref{prelWeights} we recall known facts about real representations of compact groups
used in what follows.
Section~\ref{Prel} discusses:
1) systems of random real $\pi$-polynomials
and the expectation of the number of their roots;
2) the asymptotics of the expectation for a growing representation $\pi$;
3) the geometry of ellipsoids arising in these questions; and
4) some special properties of ellipsoids in the case of a simple group $K$.
In Section~\ref{proportion}
we (1) formulate and prove the BKK theorem for groups in the form required here,
(2) use it to compute the probability of a real root of a system of random $\pi$-polynomials,
and (3) use this computation to determine the limiting probability of a real root in Theorem~\ref{thmAntiKac2}.
\subsection{Facts on representations of compact groups}\label{prelWeights}
In what follows, we use the following notions and facts from group theory:

$\bullet$\ $K$, $\mathfrak k$, and $\mathfrak k^*$ denote a connected compact Lie group,
its Lie algebra, and the space of linear functionals on $\mathfrak k$, respectively;

$\bullet$\ $T^k$, $\mathfrak t$, and $\mathfrak t^*$ denote a maximal torus in $K$,
its Lie algebra, and the space of linear functionals on $\mathfrak t$, respectively;

$\bullet$\ $\Z^k$ is the lattice of characters in $\mathfrak t^*$,
i.e., the lattice consisting of differentials of characters of the torus $T^k$;

$\bullet$\ $W^*$ is the Weyl group acting in the space $\mathfrak t^*$,
and $|W|$ denotes the number of its elements;

$\bullet$\ $\mathfrak C^*$ is a Weyl chamber in $\mathfrak t^*$;

$\bullet$\ $R$ is the root system in $\mathfrak t$, and $R^+$ is the set of positive roots;

$\bullet$\ $\rho=\frac{1}{2}\sum_{\beta\in R^+}\beta$ is the half-sum of positive roots;

$\bullet$\ $\tau=(\cdot,\cdot)$ denotes the invariant scalar product on $\mathfrak k$,
as well as the dual scalar product on $\mathfrak k^*$ and their
restrictions to $\mathfrak t$ and $\mathfrak t^*$;
$\vol_\tau(K)$ and $\vol_\tau(T^k)$ are the corresponding volumes of $K$ and $T^k$;

$\bullet$\ $d\nu$ denotes the Lebesgue measures in the spaces $\mathfrak k$,
$\mathfrak k^*$, and $\mathfrak t^*$ corresponding to the metric $\tau$;

$\bullet$\ $P(\lambda)=\prod_{\beta\in R^+}(\lambda,\beta)$;

$\bullet$\ $\mu_\lambda$ denotes the irreducible representation with highest weight $\lambda$;

$\bullet$\ $K^\C$ denotes the complexification of the group $K$;

$\bullet$\ $\mathfrak k_\C=\mathfrak k+i\mathfrak k$ is the complexification of the Lie algebra $\mathfrak k$.
\par\smallskip
The following statement follows from standard properties of the group $W^*$.
\begin{assert}\label{prPair}
For any $\lambda\in\mathfrak C^*$ there exists a unique $\lambda'\in\mathfrak C^*$
such that $\lambda'\in W^*(-\lambda)$.
\end{assert}
\noindent
Note that $\lambda''=\lambda$.
If $\lambda\in W^*(-\lambda)$, then $\lambda'=\lambda$.
For example, if the group $W^*$ contains the central symmetry map, then $\lambda'=\lambda$.
\begin{definition}\label{dfSymm}
An unordered pair $(\lambda,\lambda')$ is called a \emph{symmetric pair}.
A subset $\Lambda\subset\mathfrak C^*$ is called \emph{symmetric}
if $\lambda\in\Lambda$ implies $\lambda'\in\Lambda$.
For a symmetric set $\Lambda$ denote by $\Lambda'$
the set of symmetric pairs $\{(\lambda,\lambda'):\lambda\in\Lambda\}$.
\end{definition}
For example,
if a set $B\subset\mathfrak k^*$ is invariant under the action of $W^*$ and centrally symmetric,
then $B\cap\mathfrak C^*$ is symmetric.
For $K=T^n$, the symmetry condition for a pair $(\lambda,\delta)$ is simply $\delta=-\lambda$,
and the symmetry condition for a set means
that the set is centrally symmetric.

The following statement is an analogue of the highest-weight theory for real representations of the group $K$;
see \cite[Chapter IX, Appendix II]{B}.
\begin{assert}\label{prB}
There is a mapping $(\lambda,\lambda')\mapsto\pi_{\lambda,\lambda'}$
from the set of symmetric pairs $(\lambda,\lambda')$ in $(\mathfrak C^*\cap\Z^k)\times(\mathfrak C^*\cap\Z^k)$
onto the set of irreducible real representations of the group $K$.
This mapping is bijective.
It classifies irreducible representations into three types: \emph{real}, \emph{complex}, and \emph{quaternionic}:

(i)\ real type $\pi_{\lambda,\lambda'}$: $\:\pi_{\lambda,\lambda'}\otimes_\R\C=\mu_\lambda$ and $\lambda=\lambda'$;

(ii)\ quaternionic type $\pi_{\lambda,\lambda'}$: $\:\pi_{\lambda,\lambda'}\otimes_\R\C=\mu_\lambda\oplus\mu_\lambda$ and $\lambda=\lambda'$;

(iii)\ complex type $\pi_{\lambda,\lambda'}$: $\:\pi_{\lambda,\lambda'}\otimes_\R\C=\mu_\lambda\oplus\mu_{\lambda'}$ and $\lambda\ne\lambda'$.
\end{assert}
\section{Average number of real roots}
\label{Prel}
Throughout we use the following notation:

$\bullet$\
$\tau$ denotes the metric on the Lie algebra $\mathfrak k$ of the group $K$, invariant under the adjoint action of $K$,
such that the corresponding invariant measure $\chi$ on $K$ is probabilistic, i.e., $\int_K d\chi = 1$;

$\bullet$\
the notation $\tau$ is also used for the dual metric on $\mathfrak k^*$ as well as for the restrictions of the metrics on $\mathfrak k, \mathfrak k^*$ to the subspaces $\mathfrak t$ and $\mathfrak t^*$;

$\bullet$\
$\nu$ denotes the Lebesgue measure corresponding to the metric $\tau$ on the Lie algebra $\mathfrak k$,
and the corresponding Lebesgue measures on the spaces
$\mathfrak t,\mathfrak t^*,\mathfrak k^*$;
these measures are used below for measuring the volumes of convex bodies;

$\bullet$\
the scalar product $(*,*)$ in the space of $\pi$-polynomials ${\rm Trig}(\pi)$
is induced from the space
$L^2_\R(\chi)$.
\subsection{Average number of roots (Definition)}\label{s2.1}
Let $\pi$ be a finite-dimensional real representation of a connected compact Lie group $K$.
Using the Gaussian measure
$$
 \mu_\pi(U) = \frac{1}{(2\pi)^{\frac{\dim {\rm Trig}(\pi)}{2}}}
 \int_U \exp\!\big( - (f,f)/2 \big) \, df
$$
on the space ${\rm Trig}(\pi)$,
we regard $\pi$-polynomials $f_1,\ldots,f_n$ as
independent normally distributed random vectors in ${\rm Trig}(\pi)$.
Denote by $\mathfrak M(\pi)$ the expectation of the random variable
equal to the number of roots of the random system of equations $f_1 = \ldots = f_n = 0$.

An equivalent way of defining $\mathfrak M(\pi)$ is as follows.
Denote by $\P_i$ the projectivization of the space ${\rm Trig}(\pi_i)$,
i.e., the set of one-dimensional subspaces in ${\rm Trig}(\pi_i)$.
For any $\lambda_1,\ldots,\lambda_n \in \R \setminus 0$,
the roots of all systems of equations of the form $\lambda_1 f_1 = \ldots = \lambda_n f_n = 0$ coincide.
Consider the number of these roots as a function $N(f_1,\ldots,f_n)$ on the product $\P_1\times\ldots\times\P_n$
of the projective spaces $\P_i$.
Let $\phi_i$ be the orthogonally invariant probability measure on the projective space $\P_i$.
Such a measure exists and is unique.
\begin{definition}\label{dfAverage}
\begin{equation}\label{eqAverage}
\mathfrak M(\pi_1,\ldots,\pi_n)
   = \int_{(f_1,\ldots,f_n) \in \P_1 \times \ldots \times \P_n}
     N(f_1,\ldots,f_n) \: d\phi_1 \ldots d\phi_n,
\end{equation}
is called the \emph{average number of roots of systems of random $\pi_i$-polynomials}.
We also use the notation
$\mathfrak M(\pi) = \mathfrak M(\pi_1,\ldots,\pi_n)$ when
$\pi = \pi_1 = \ldots = \pi_n$.
\end{definition}
\begin{definition}\label{dfFlat}
A representation $\mu$ without multiple irreducible components is called \emph{flat}.
For any real or complex representation $\mu$,
the flat representation $\mu_F$ with the same irreducible components
is called the \emph{flattening} of $\mu$.
\end{definition}
The following two statements are direct consequences of
Definitions~\ref{dfFlat} and~\ref{dfAverage}.
\begin{corollary}\label{corFl1}
For any $\mu$ one has ${\rm Trig}(\mu) = {\rm Trig}(\mu_F)$.
\end{corollary}
\begin{corollary}\label{corFl2}
Let $\pi_{1F},\ldots,\pi_{nF}$ be the flattenings of the representations $\pi_1,\ldots,\pi_n$.
Then
$\mathfrak M(\pi_1,\ldots,\pi_n) = \mathfrak M(\pi_{1F},\ldots,\pi_{nF})$.
\end{corollary}
\subsection{Newton ellipsoids}\label{s2.2}
Define the mapping $\Theta(\pi)\colon K\to{\rm Trig}(\pi)$ by
\begin{equation}\label{eqTheta}
\forall f\in{\rm Trig}(\pi)\colon\:\left(\Theta(\pi)(\rho),f\right)=\frac{1}{\sqrt N}f(\rho),
\end{equation}
where $N=\dim{\rm Trig}(\pi)$.
\begin{lemma}\label{lmMain1}
The set $\Theta(\pi)(K)$ lies on the sphere $S$ in the space ${\rm Trig}(\pi)$ of radius $1$ centered at the origin.
\end{lemma}
\begin{proof}
The inner product $(*,*)$ in ${\rm Trig}(\pi)$ is invariant under the action of the group $K$.
Hence the set $\Theta(\pi)(K)$ is contained in a sphere of some radius $r$.
From (\ref{eqTheta}) it follows that for any orthonormal basis
 $f_1,\ldots,f_N$
 in the space ${\rm Trig}(\pi)$ we have
\begin{equation}\label{eqTheta2}
  \Theta(\pi)(\rho)= \frac{1}{\sqrt N}(f_1(\rho)f_1+\ldots+f_N(\rho)f_N).
\end{equation}
Therefore
$
r^2=\left(\Theta(\pi)(\rho),\Theta(\pi)(\rho)\right)=\frac{1}{N}\sum_if^2_i(\rho).
$
Integrating the function $\sum_if^2_i(\rho)$
with respect to the measure $d\chi$, we get
$$
r^2=\int_K r^2\:d\chi=\frac{1}{N}\sum_i\int f^2_i(\rho)\:d\chi=\left((f_1,f_1)+\ldots+(f_N,f_N)\right)/N=1.
$$
\end{proof}
\begin{definition}\label{dfKill_pi}
Define a symmetric bilinear form $\mathcal G(\pi)$ in the tangent space $T_\rho K$ of the group $K$ at the point $\rho$
as the pullback of the inner product $(*,*)$ in ${\rm Trig}(\pi)$
under the mapping $\Theta(\pi)$.
Let $F(\pi)$ be the restriction of $\mathcal G(\pi)$ to the Lie algebra $\mathfrak k$,
and let $g(\pi)$ be the corresponding quadratic form on $\mathfrak k$.
\end{definition}

From equality (\ref{eqTheta2}) the following statement follows.
\begin{corollary}\label{corFbasis}
For any orthonormal basis $f_1,\ldots,f_N$ in ${\rm Trig}(\pi)$,
$$F(\pi)(\xi,\eta)=\frac{1}{N}\sum_{1\leq i\leq N} df_i(\xi)\: df_i(\eta),$$
where $df_i$ is the differential of the function $f_i$ at the identity element of the group $K$.
\end{corollary}
The quadratic form $g(\pi)$ is non-negative.
Set $h_\pi(x)=\sqrt{g(\pi)(x)}$.
The function $h_\pi$ is non-negative,
positively homogeneous of degree $1$, and convex.
Hence $h_\pi$ is the support function of a uniquely determined compact convex set
${\rm Ell}(\pi)$ in $\mathfrak k^*$,
i.e., $h_\pi(x)=\max_{y \in{\rm Ell}(\pi)} y(x)$.
\begin{definition}\label{dfSuppRepr}
The set ${\rm Ell}(\pi)\subset\mathfrak k^*$ is called the Newton ellipsoid of the representation $\pi$.
It is centrally symmetric and
is an ellipsoid in the subspace orthogonal to the kernel of the quadratic form $g(\pi)$.
\end{definition}
The following statement follows from the invariance of $h_\pi$ with respect to the adjoint action of $K$.
\begin{corollary}\label{cor21-1}
The ellipsoid ${\rm Ell}(\pi)$ is invariant under the coadjoint action of $K$ on $\mathfrak k^*$.
\end{corollary}
\begin{corollary}\label{Ball}

{\rm (1)}
The ellipsoid ${\rm Ell}(\pi)$ is a ball of some
invariant metric in the space ${\rm ker}^\bot\subset\mathfrak k^*$,
where ${\rm ker}^\bot$ is the orthogonal complement to the subspace ${\rm ker}(d\pi)\subset\mathfrak k$.

{\rm (2)} If the group $K$ is simple,
then the ellipsoid ${\rm Ell}(\pi)$ is a ball in $\mathfrak k^*$ for any coadjointly invariant metric in $\mathfrak k^*$.
\end{corollary}
\begin{proof}
Statement (1) follows from Corollary \ref{cor21-1}.
If the group $K$ is simple, then any two invariant metrics differ by a constant factor;
this implies (2).
\end{proof}
\begin{corollary}\label{cor21-2}
For the flattening $\pi_F$ of the representation $\pi$, we have ${\rm Ell}(\pi_F)={\rm Ell}(\pi)$.
\end{corollary}
\begin{proof}
Follows from Corollary \ref{corFl1}.
\end{proof}
\begin{example}\label{exAnti2}
Consider the representation $\mu$ from Example \ref{exKac}.
For $k=1,\ldots,m$,
the basis $1$, $\sqrt 2\cos(k\theta)$,
 $\sqrt 2\sin(k\theta)$
of the space ${\rm Trig}(\mu)$ is orthonormal.
Since $\frac{d}{d\theta}\cos \theta=-2\pi i\sin\theta$,
$\frac{d}{d\theta}\sin \theta=2\pi i\cos\theta$,
it follows from Corollary \ref{corFbasis} and Definition \ref{dfSuppRepr}
that the support function of the Newton ellipsoid ${\rm Ell}(\mu)$ is
$$
 h_\mu(\xi)=\sqrt{g_\pi(\xi)}=2\pi \vert\xi\vert\sqrt{2 (1^2+\ldots+m^2)}/\sqrt{2m+1}=2\pi\vert\xi\vert\sqrt{m(m+1)/3}.
$$
The Newton ellipsoid ${\rm Ell}(\mu)$ is the interval with endpoints at
$\pm h_\mu(1)$.
Since $h_\mu(1)=2\pi\sqrt{m(m+1)/3}$,
it follows from Theorem \ref{thmMixed} below that
$\mathfrak M(\pi)=2\sqrt{m(m+1)/3}$.
\end{example}
\subsection{Mixed volume of Newton ellipsoids}\label{s2.3}
\begin{theorem}\label{thmMixed}
For finite-dimensional real representations $\pi_1,\ldots,\pi_n$ of the group $K$,
the following holds:
$$
\mathfrak M(\pi_1,\ldots,\pi_n)=\frac{n!}{(2\pi)^n}\:\vol_\tau({\rm Ell}(\pi_1),\ldots,{\rm Ell}(\pi_n)).
$$
\end{theorem}
\noindent
\emph{Proof.}
Recall that a convex body $\mathcal B$ on a smooth manifold $X$
is defined as a family of centrally symmetric compact convex subsets $\mathcal B(x)$
in the fibers $T^*_xX$ of the cotangent bundle of $X$;
see \cite{AK, K20}.
Define the volume $\vol(\mathcal B)$ of the convex body $\mathcal B$ as
the symplectic volume of the set $\bigcup_{x\in X}\mathcal B(x)\subset T^*X$.
More precisely, the volume is measured using the volume form $\omega^n/n!$,
where $\omega$ is the standard symplectic form on the cotangent bundle
$T^*X$;
see \cite{Ar}.

For compact convex sets in a vector space,
the operations "addition" (Minkowski addition) and "multiplication by a non-negative number $\lambda$" (homothety with coefficient $\lambda$) are defined.
The Minkowski sum of convex sets $A$ and $B$ is the set formed by the pairwise sums of points $A$ and $B$.
Recall that the Minkowski sum of convex sets $A$ and $B$ is the set
formed by all pairwise sums of points from $A$ and $B$.
The Minkowski addition satisfies the cancellation law.
Using Minkowski sums and homotheties in the fibers $T^*_xX$,
we define, for nonnegative $\lambda_i$,
the linear combination $\lambda_1\mathcal B_1+\ldots+\lambda_k\mathcal B_k$
of convex bodies on $X$ as
$$
(\lambda_1\mathcal B_1+\ldots+\lambda_k\mathcal B_k)(x)=\lambda_1\mathcal B_1(x)+\ldots+\lambda_k\mathcal B_k(x).
$$
\begin{definition}\label{dfMixed}
If the volumes of convex bodies
$\mathcal B_1,\ldots,\mathcal B_k$
on $X$
are finite,
then $\vol(\lambda_1{\mathcal B}_1+ \ldots+ \lambda_k{\mathcal B}_k)$ is
a homogeneous polynomial of degree $n$ in the variables $\lambda _1, \ldots, \lambda _k$,
where $n=\dim X$.
The coefficient
of the volume polynomial
$\vol(\lambda_1{\mathcal B}_1+ \ldots+ \lambda_n{\mathcal B}_n)$
at the monomial $\lambda _1\cdot\ldots\cdot \lambda _n$,
divided by $n!$,
is denoted by $\mathcal V({\mathcal B}_1,\ldots,\mathcal B_n)$ and called the \emph{mixed volume} of $n$ convex bodies on $X$.
\end{definition}
Let $V\subset C^\infty(X)$ be a finite-dimensional vector space
such that
\begin{equation}\label{eqNot0}
  \forall x\in X\:\exists f\in V\colon\: f(x)\neq0.
\end{equation}
Assume that a scalar product $(*,*)$ is chosen in the space $V$.
Consider the map $\theta\colon X\to V$ such that
$$
\forall (f\in V,\:x\in X)\colon\:(\theta(x),f)=f(x).
$$
From condition (\ref{eqNot0}) it follows that $\forall x\in X\colon\theta(x)\ne0$.
Set $\Theta(x)=\theta/\vert\theta(x)\vert$.
\begin{definition}\label{defAK}
Define the convex body $\mathcal B_V$ on $X$ by $\mathcal B_V(x)=d_x^*\Theta(B)$,
where $B\subset V$ is the ball of radius $1$ centered at the origin,
and $d_x^*\Theta\colon V\to T^*_xV$ is the linear map adjoint to the differential $d\Theta$
at the point $x\in X$.
For each $x\in X$, the set $\mathcal B_V(x)$ is an ellipsoid at $x$.
\end{definition}
Let $V_1,\ldots,V_n\subset C^\infty(X)$ be finite-dimensional vector spaces with fixed scalar products in $C^\infty(X)$.
As at the beginning of \textsection\ref{s2.1},
the average number of common zeros
$\mathfrak M(V_1,\ldots,V_n)$
of random functions $f_1\in V_1,\ldots,f_n\in V_n$ is defined; see \cite{AK}.
In \cite[Theorem 1]{AK} it is proved that,
if for each space $V_i$ condition (\ref{eqNot0}) holds,
then
\begin{equation}\label{eqAK}
  \mathfrak M(V_1,\ldots,V_n)=\frac{n!}{(2\pi)^n}\mathcal V(\mathcal B_{V_1},\ldots,\mathcal B_{V_n}).
\end{equation}

Let $X=K$ and $V_i={\rm Trig}(\pi_i)$.
The Newton ellipsoid ${\rm Ell}(\pi)$ is located in the cotangent space $T^*_eK$ at the identity of $K$.
Denote by $\mathcal B(\pi)$ the convex body on $K$
consisting of left shifts of the ellipsoid ${\rm Ell}(\pi)$.
From the definition of the Newton ellipsoid and from Lemma \ref{lmMain1} it follows that $\mathcal B(\pi_i)=\mathcal B_{V_i}$.
Applying (\ref{eqAK}) we obtain
$$
\mathfrak M(\pi_1,\ldots,\pi_n)=\frac{n!}{(2\pi)^n}\mathcal V(\mathcal B(\pi_1),\ldots,\mathcal B(\pi_n)).
$$
By the left-invariance of the convex body $\mathcal B(\pi_i)$ on $K$
it follows that
$$
\mathfrak M(\pi_1,\ldots,\pi_n)=\frac{n!}{(2\pi)^n}\vol_\tau(\mathcal B(\pi_1),\ldots,\mathcal B(\pi_n)).
$$
The theorem is proved.
\par\smallskip

The following two statements are not used later.
They are included for completeness.
\begin{corollary}\label{corEqui}
Let $U$ be an open subset of $K$.
Then the average number of common zeros of $n$ random $\pi_i$-polynomials contained in $U$
is equal to $$\frac{n!}{(2\pi)^n}\:\vol_\tau({\rm Ell}(\pi_1),\ldots,{\rm Ell}(\pi_n))\:\int_U d\chi.$$
\end{corollary}
\begin{proof}
Equality (\ref{eqAK}) from \cite[Theorem 1]{AK} remains valid when restricting functions from the spaces $V_i$ and convex bodies $\mathcal B_i$
to any open subset $U$ in $X$.
The statement follows.
\end{proof}
\begin{corollary}\label{Hodge}
For any representations $\pi_1,\ldots,\pi_n$
\begin{itemize}
\item[(1)]
$
  \mathfrak M^2(\pi_1,\ldots,\pi_n)\geq\mathfrak M(\pi_1,\ldots,\pi_{n-1},\pi_{n-1})\cdot\mathfrak M(\pi_1,\ldots,\pi_n,\pi_n),
$
\item[(2)] ${\mathfrak M}^n(\pi_1,\ldots,\pi_n) \ge \mathfrak M(\pi_1)\cdot\ldots\cdot \mathfrak M(\pi_n)$.
\item[(3)]
If the group $K$ is simple,
then all of these inequalities become equalities.
\end{itemize}
\end{corollary}
\begin{proof}
(1) and (2) follow from the Alexandrov–Fenchel inequalities for the mixed volumes
of the ellipsoids ${\rm Ell}(\pi_i)$;
see \cite{Al}.
The mixed volume of balls equals the product of their radii,
multiplied by the volume of the unit ball.
Hence (3) follows from Corollary \ref{Ball} (2).
\end{proof}
\begin{remark}
The inequalities of Corollary \ref{Hodge} are analogues of the Hodge inequalities for intersection indices of hypersurfaces
in a projective algebraic variety;
see, e.g., \cite{KK}.
\end{remark}
\subsection{Using complexification}\label{sMean1}
Consider the complex vector space of $\mu$-polynomials ${\rm Trig}(\mu)$,
where $\mu$ is a complex representation of the group $K$.
Assume that ${\rm Trig}(\mu)$ is equipped with a Hermitian metric,
inherited from the space $L^2_\C(d\chi)$,
where $\chi$ is the invariant measure on $K$.
\begin{lemma}\label{lmReIm}
Let $\pi$ be a real representation of $K$, and $\mu=\pi\otimes_\R\C$.
Then

{\rm(1)}\ ${\rm Trig}(\pi)$ is a real subspace of ${\rm Trig}(\mu)$;

{\rm(2)}\ The restriction of the Hermitian inner product $\langle*,*\rangle$ to ${\rm Trig}(\pi)$
coincides with the scalar product taken from $L^2_\R(d\chi)$;

{\rm(3)}\ An orthonormal basis of ${\rm Trig}(\pi)$ is also an orthonormal basis of ${\rm Trig}(\mu)$;

{\rm(4)}\ Denote by $\re(f)$ the real part of a function $f\colon K\to\C$,
and set $\re{\rm Trig}(\mu)=\{\re(f)\colon f\in{\rm Trig}(\mu)\}$.
Then ${\rm Trig}(\pi)=\re{\rm Trig}(\mu)$.
\end{lemma}
\noindent
\emph{Proof.}
All four statements are direct consequences of the definition of the complexification $\pi\otimes_\R\C$ of the representation $\pi$
and of the Hermitian inner product $\langle *,*\rangle$ in the space ${\rm Trig}(\mu)$.
\begin{corollary}\label{corRP2}
For the representation $\pi_{\lambda,\lambda'}$ (see Proposition \ref{prB}) we have
${\rm Trig}(\pi_{\lambda,\lambda'})=\re{\rm Trig}(\mu_\lambda)$.
\end{corollary}
\begin{proof}
From the definition of duality of representations it follows that
for a pair of dual complex representations $\mu$ and $\nu$,
$\re{\rm Trig}(\mu)=\re{\rm Trig}(\nu)$.
Consequently (see \textsection\ref{s1.2}) $\re{\rm Trig}(\mu_\lambda)=\re{\rm Trig}(\mu_{\lambda'})$.
The statement now follows from Lemma \ref{lmReIm} (4) and Proposition \ref{prB}.
\end{proof}
For a complex representation $\mu$, similarly to the real case in (\ref{eqTheta}),
define the mapping $\Theta(\mu)\colon K\to{\rm Trig}(\mu)$
\begin{equation}\label{eqThetaC}
\forall f\in{\rm Trig}(\mu)\colon\:\langle f,\Theta(\mu)(g)\rangle=\frac{1}{\sqrt N}f(g),
\end{equation}
where $N=\dim_\C{\rm Trig}(\mu)$.
As in the real case (see (\ref{eqTheta2})),
for any orthonormal basis $f_1,\ldots,f_N$ of ${\rm Trig}(\mu)$, we have
\begin{equation}\label{eqThetaC2}
  \Theta(\mu)(g)= \frac{1}{\sqrt N}(f_1(g)f_1+\ldots+f_N(g)f_N).
\end{equation}
Similarly to Definition \ref{dfKill_pi},
consider the complex-valued bilinear form
in the space $\mathfrak k$
\begin{equation}\label{eqKill}
 F(\mu)(\xi,\eta)=\langle d\Theta(\mu)(\xi),d\Theta(\mu)(\eta)\rangle=\frac{1}{N}\sum_i df_i(\xi)\:\overline{df_i(\eta)}.
\end{equation}
The form $F(\mu)$ is Hermitian symmetric, i.e., $F(\mu)(\xi,\eta)=\overline{F(\mu)(\eta,\xi)}$.
From Lemma \ref{lmReIm} (3) it follows that
\begin{equation}\label{eqThetaC3}
  F(\pi\otimes_\R\C)=F(\pi).
\end{equation}
Similarly to the real case,
\begin{equation}\label{eqFlComplex}
  F(\mu_F)=F(\mu),
\end{equation}
where $\mu_F$ is the flattening of the representation $\mu$;
see Definition \ref{dfFlat}.

Below we use the concepts from \textsection\ref{prelWeights}.
Let $\Lambda$ be a finite symmetric subset of $\Z^k\cap\mathfrak C^*$,
and
\begin{equation}\label{pimu}
  \pi=\bigoplus_{(\lambda,\lambda')\in\Lambda'}m_\lambda\pi_{\lambda,\lambda'},
\,\,\,\,\mu=\bigoplus_{\lambda\in\Lambda}m_\lambda\mu_\lambda.
\end{equation}
Recall that $\pi_{\lambda,\lambda'}$ and $\mu_\lambda$ are irreducible real and complex representations
of the group $K$, respectively.
The sets $\Lambda$ and $\Lambda'$ are called the spectra,
and their elements the weights,
of the representations $\mu$ and $\pi$.
From Proposition \ref{prB} it follows that
$$
\pi\otimes_\R\C=\mu+\sum_{\lambda\in \Q(\Lambda)}m_\lambda\mu_\lambda,
$$
where $\Q(\Lambda)$ is the set of $\lambda\in\Lambda$
for which the weight $(\lambda,\lambda)$ is quaternionic.
Hence
$\left(\pi\otimes_\R\C\right)_F=\mu_F$.
From (\ref{eqFlComplex}) and (\ref{eqThetaC3}) it follows that
\begin{equation}\label{eqFlPi}
F(\pi)=F(\mu).
\end{equation}
Below in proofs
we assume,
based on (\ref{eqFlComplex}) and on Corollary \ref{cor21-2},
that all representations are flat.
\begin{lemma}\label{lmFhigh}
For the representation $\mu_\lambda$ with highest weight $\lambda$ we have
$$
F(\mu_\lambda)(\xi,\eta)=-\frac{1}{\dim\mu_\lambda}\tr\left(d\mu_\lambda(\xi)\cdot d\mu_\lambda(\eta)\right),
$$
where $d\mu_\lambda$ is regarded as a representation of the Lie algebra $\mathfrak k$.
\end{lemma}
\begin{proof}
Represent the operators of $\mu_\lambda$
by unitary matrices $\{t^\lambda_{i,j}\}$.
From the orthogonality relations for the matrix elements $t^\lambda_{i,j}$
(see, e.g., \cite{B})
it follows that
the functions $\sqrt {\dim(\mu_\lambda)}\: t^\lambda_{i,j}$
form an orthonormal basis of the space
${\rm Trig}(\mu_\lambda)$.
According to (\ref{eqKill}),
$$
F(\mu_\lambda)(\xi,\eta)=\frac{1}{\dim(\mu_\lambda)}\sum_{i,j} dt^\lambda_{i,j}(\xi)\overline {dt^\lambda_{i,j}(\eta)}=-\frac{1}{\dim\mu_\lambda}\tr\left(d\mu_\lambda(\xi)\cdot d\mu_\lambda(\eta)\right).
$$
The lemma is proved.
\end{proof}
The adjoint representation of a simple Lie group $K$ is irreducible.
Therefore, it is equal to $\mu_\alpha$,
where $\alpha$ is an element of the Weyl chamber $\mathfrak C^*$, called the highest root of the group $K$.
Recall that the symmetric bilinear form $\kappa(\xi,\eta)=-\tr\left(d\mu_\alpha(\xi)\cdot d\mu_\alpha(\eta)\right)$ on the Lie algebra $\mathfrak k$
is non-degenerate, positive, and is called the Killing metric.
\begin{corollary} \label{corKilling}
Let $\alpha$ be the highest root of a simple Lie group $K$ of dimension $n$.
Then (for the element $\alpha'\in\mathfrak C^*$ see Definition \ref{dfSymm})
$$\forall \xi,\eta\in\mathfrak k\colon\,F(\pi_{\alpha,\alpha'})(\xi,\eta)=F(\mu_\alpha)(\xi,\eta)=n^{-1}\kappa(\xi,\eta).$$
\end{corollary}
\begin{lemma}\label{lmMulty}
Let $p(\lambda)=\dim(\mu_\lambda)$. Then
$$
F(\pi)=\frac{\sum_{\lambda\in\Lambda} p^2(\lambda)\:F(\mu_\lambda)}{\sum_{\lambda\in\Lambda} p^2(\lambda)}.
$$
\end{lemma}
\noindent
\emph {Proof.}
Since $\dim {\rm Trig}(\mu_\lambda)=p^2(\lambda)$, we have $\dim\left({\rm Trig}(\mu)\right)=\sum_{\lambda\in\Lambda} p^2(\lambda)$.
From  (\ref{eqKill}) it follows that
$
F(\mu)=\left(\sum_{\lambda\in\Lambda} p^2(\lambda)\:F(\mu_\lambda)\right)/\sum_{\lambda\in\Lambda} p^2(\lambda).
$
The statement now follows from equality (\ref{eqFlPi}).
\begin{corollary}\label{corUnitary}
For any $\xi,\eta\in\mathfrak k$,
$$
F(\pi)(\xi,\eta)=-\frac{\sum_{\lambda\in\Lambda} p(\lambda)\:\tr(d\mu_\lambda(\xi)\cdot d\mu_\lambda(\eta))}{\sum_{\lambda\in\Lambda} p^2(\lambda)}.
$$
\end{corollary}
\begin{proof}
This follows from Lemmas  \ref{lmMulty}, \ref{lmFhigh} and from equality (\ref{eqFlPi}).
\end{proof}
\subsection{Simple groups}\label{sMean2}
For a simple group $K$, the Newton ellipsoid ${\rm Ell}(\pi)$ is a ball centered at the origin;
see Corollary \ref{Ball}, (2).
Here we compute the radius of this ball.
\begin{theorem}\label{thmRadius}
Let $r(\Lambda)$ be a radius of the ball ${\rm Ell}(\pi)$.
Then
$$
r^2(\Lambda)=\frac{\sum_{\lambda\in\Lambda} p^2(\lambda)\:(\lambda,\lambda+2\rho)}{n(\alpha,\alpha+2\rho)\sum_{\lambda\in\Lambda} p^2(\lambda)},
$$
where $\alpha$ is the highest root of the group $K$,
i.e., the highest weight of the adjoint representation $\mu_\alpha$.
\end{theorem}
For a simple group, any two invariant metrics differ by a constant factor.
Hence, for any complex representation $\mu$, there exists a number $l(\mu)>0$
such that for any $\xi,\eta\in\mathfrak k$
\begin{equation}\label{ind1}
\frac{\tr\left(d\mu(\xi)\cdot d\mu(\eta)\right)}{\dim(\mu)}=l(\mu)\frac{\tr\left(d\mu_\alpha(\xi)\cdot d\mu_\alpha(\eta)\right)}{n}.
\end{equation}
It is known (see, e.g., \cite[(8)]{Vinb}) that
\begin{equation}\label{ind2}
l(\mu_\lambda)=\frac{(\lambda,\lambda+2\rho)}{(\alpha,\alpha+2\rho)},
\end{equation}
where $\rho$ is the half-sum of positive roots of the group $K$; see \textsection\ref{prelWeights}.
\begin{remark}
Equality (\ref{ind2}) is equivalent to the fact that $(\lambda,\lambda+2\rho)$
is the eigenvalue of the Casimir operator in the space ${\rm Trig}(\mu_\lambda)$; see \cite{B}.
\end{remark}
\begin{lemma}\label{lmIndex}
Let $\pi$ be a representation of the group $K$ with spectrum $\Lambda'$
defined by (\ref{pimu}).
Then
$$
F(\pi)(\zeta,\zeta)=
\frac{|\zeta|^2}{n(\alpha,\alpha+2\rho)}\:
\frac{\sum_{\lambda\in\Lambda} p^2(\lambda)\:(\lambda,\lambda+2\rho)}{\sum_{\lambda\in\Lambda} p^2(\lambda)}.
$$
\end{lemma}
\begin{proof}
Applying first Corollary \ref{corUnitary} and then
equalities (\ref{ind1}), (\ref{ind2}), we obtain
\begin{equation*}
 F(\pi)(\zeta,\zeta)=-\frac{\sum p(\lambda)\:\tr(d\mu_\lambda(\zeta)\cdot d\mu_\lambda(\zeta))}{\sum p^2(\lambda)}
  =
\frac{\vert\zeta\vert^2}{n(\alpha,\alpha+2\rho)}\frac{\sum p^2(\lambda)\:(\lambda,\lambda+2\rho)}{\sum p^2(\lambda)}.
\end{equation*}
\end{proof}
\par\smallskip
\noindent
\emph{Proof of Theorem \ref{thmRadius}.}
Recall that $h_\pi=\sqrt{F(\pi)}$ is the support function of the Newton ellipsoid ${\rm Ell}(\pi)$;
see Definition \ref{dfSuppRepr}.
Since ${\rm Ell}(\pi)$ is a ball,
its radius $r(\Lambda)$ equals $h_\pi(\zeta)$ for any $\vert\zeta\vert=1$.
Therefore, the theorem follows from Lemma \ref{lmIndex}.
\subsection{Asymptotics of the average number of roots}\label{sAsymp}
Here we also assume that the group $K$ is simple.
Let $\Delta$ be a compact convex set in $\mathfrak t^*$.
Assume that $\Delta$ is centrally symmetric and invariant under the action of the Weyl group.
Examples of such sets are balls centered at the origin or weight polytopes of real representations of the group $K$;
see Definition \ref{dfweightedpolyhedron}.
The finite set $\Lambda=\Delta\cap\mathfrak C^*\cap\Z^k$ is symmetric in the sense of Definition \ref{dfSymm}.
As $m\to\infty$,
we consider the representation
$$\pi_m=\bigoplus_{(\lambda,\lambda')\in (m\Lambda')}\pi_{\lambda,\lambda'},$$
(which coincides with (\ref{eqpim}) for $\Lambda=B$), and
compute the asymptotics of the average number of roots $\mathfrak M(\pi_m)$
as $m\to\infty$ (Theorem \ref{thmAsymp1}). We will later use this computation
to prove Theorem \ref{thmAntiKac2}.

Recall that the Newton ellipsoids
${\rm Ell}(\pi_m)$ are balls in $\mathfrak k^*$ centered at the origin;
see Corollary \ref{Ball} (2).
The following statement computes the asymptotic growth of the radius
of the ball
${\rm Ell}(\pi_m)$ as $m\to\infty$.
\begin{theorem}\label{thmAsymp1}
If the group $K$ is simple,
and $\Delta$ is the unit ball in $\mathfrak t^*$,
then
$$
\lim_{m\to\infty}\frac{\mathfrak M(\pi_m)}{m^n}=\frac{n!}{(2\pi)^n}\sigma_n (n+2)^{-n/2}(\alpha,\alpha+2\rho)^{-n/2},
$$
where $n=\dim K$, $\sigma_n$ is the volume of the $n$-dimensional unit ball,
$\alpha$ is the highest root, and $\rho$ is the half-sum of the positive roots of the group $K$.
\end{theorem}
Denote by $\mathfrak N(\Delta)$
the compact set in $\mathfrak k^*$ consisting of coadjoint orbits of $K$ intersecting $\Delta$.
From the results of \cite{K},
it follows that the set $\mathfrak N(\Delta)$ is convex.
In particular, if $\Delta$ is a ball in the space $\mathfrak t^*$ centered at the origin,
then $\mathfrak N(\Delta)$ is a ball of the same radius in the space $\mathfrak k^*$.
More details on the properties of $\mathfrak N(\Delta)$ can be found in \textsection\ref{BKK}.
\begin{theorem}\label{thmAsymp}
As $m\to\infty$, the sequence of balls $\frac{1}{m}{\rm Ell}(\pi_m)$
converges in the Hausdorff topology
to a ball of radius $r_\Delta$,
where
\begin{equation}\label{eqAsymp1}
  r^2_\Delta=\frac{\int_{\mathfrak N(\Delta)}(\xi,\xi)}
{n(\alpha,\alpha+2\rho)\:\vol(\mathfrak N(\Delta))}
\end{equation}
\end{theorem}
\begin{proof}
It is known that the dimension $p(\lambda)$ of the irreducible representation $\mu_\lambda$ with highest weight $\lambda$ is equal to $\frac{P(\lambda+\rho)}{P(\rho)},$
where $P(\lambda)=\prod\nolimits_{\beta\in R^+}(\lambda,\beta)$; see \cite[Theorem 5 in Ch. IX, \S 7, n. 3]{B}).
Set $A_m=\Delta\cap\mathfrak C^*\cap\left(\frac{1}{m}\Z^k\right)$.
From Lemma \ref{lmIndex} it follows that
\begin{multline*}
\frac{1}{m^2}F(\pi_m)(\zeta,\zeta)=
\frac{1}{m^2}\frac{\vert\zeta\vert^2}{n(\alpha,\alpha+2\rho)}\frac{\sum_{\lambda\in m\Lambda} p^2(\lambda)\:(\lambda,\lambda+2\rho)}{\sum_{\lambda\in m\Lambda} p^2(\lambda)}\\
=\frac{\vert\zeta\vert^2}{n(\alpha,\alpha+2\rho)}
\frac{\sum\nolimits_{\lambda\in A_m}P^2(\lambda+\rho/m)(\lambda,\lambda+2\frac{\rho}{m})
}{\sum\nolimits_{\lambda\in A_m}P^2(\lambda+\rho/m)}
 \end{multline*}
 Hence, for large $m$, and fixed $\zeta$,
$$
\frac{1}{m^2}F(\pi_m)(\zeta,\zeta)\asymp
\frac{\vert\zeta\vert^2}{n(\alpha,\alpha+2\rho)}
\frac{\sum_{\lambda\in A_m} P^2(\lambda)
\:(\lambda,\lambda)}{\sum_{\lambda\in A_m} P^2(\lambda)}
$$
Denote by $s$ the $\tau$-volume of the fundamental parallelepiped of the character lattice $\Z^k\subset\mathfrak t^*$.
Then the sums
$$
\sum_{\lambda\in A_m} P^2(\lambda)\cdot\frac{s}{m^k},\,\,\,\,\,\,\,\sum_{\lambda\in A_m} P^2(\lambda)\:(\lambda,\lambda)\cdot\frac{s}{m^k}
$$
are Riemann sums for the integrals
$$\int_\Delta P^2(\lambda)\:d\nu,\,\,\,\,\,\,\,\int_\Delta P^2(\lambda)\:(\lambda,\lambda)\:d\nu$$
with respect to the partition $\frac{1}{m}\Z^k$.
Therefore
\begin{equation}\label{eqForWeyl}
  \lim_{m\to\infty}\frac{F(\pi_m)(\zeta,\zeta)}{m^2}=
\frac{\vert\zeta\vert^2}{n(\alpha,\alpha+2\rho)}
\frac{\int_\Delta P^2(\lambda)(\lambda,\lambda)\:d\nu}
{\int_\Delta P^2(\lambda)\:d\nu}
\end{equation}
Next we use Weyl's integration formula for a coadjoint-invariant function $f\colon\mathfrak k^*\to\R$;
see \cite[Proposition 3, Ch. IX, \S6, n. 3]{B}.
\begin{equation}\label{WeyIntFormula}
 \int_{\mathfrak k^*}f \:d\nu=\frac{\vol(K)}{\vert W\vert\vol(T^k)}\int_{\mathfrak t^*}P^2(\lambda)\:f\:d\nu,
\end{equation}
where $\vol(K),\vol(T^k)$ are the $\tau$-volumes of the group $K$ and the torus $T^k$.
Let $\varphi\colon\mathfrak k^*\to\R$ be a coadjoint-invariant function
whose restriction to $\mathfrak t^*$ is the characteristic function of the set $\Delta$.
Applying Weyl's formula to the functions $\varphi(\lambda)\:(\lambda,\lambda),\,\varphi$ from the numerator and denominator
of (\ref{eqForWeyl}),
we obtain
$$
\lim_{m\to\infty}\frac{F(\pi_m)(\zeta,\zeta)}{m^2}=\frac{\vert\zeta\vert^2\int_{\mathfrak N(\Delta)}(\xi,\xi)\:d\nu(\xi)}
{n(\alpha,\alpha+2\rho)\:\vol(\mathfrak N(\Delta))}
$$
Since $h_{\pi_m}=\sqrt {F(\pi_m)}$, we have
$$
\lim_{m\to\infty}\frac{1}{m}h_{\pi_m}(\zeta)=\sqrt{\frac{1}{n(\alpha,\alpha+2\rho)}\frac{\int_{\mathfrak N(\Delta)}(\xi,\xi)}{\vol(\mathfrak N(\Delta))}}\:\vert\zeta\vert
$$
For $\vert\zeta\vert=1$ we obtain the desired equality.
Theorem \ref{thmAsymp} is proved.
\end{proof}
\begin{corollary}\label{corAsymp1}
Let $\Delta$ be a ball of radius $r$ in $\mathfrak t^*$.
Then
\begin{equation}\label{volAsymp}
  \lim_{m\to\infty}\frac{\vol({\rm Ell}(\pi_m))}{m^n}= \sigma_n (n+2)^{-n/2}(\alpha,\alpha+2\rho)^{-n/2}\:r^n,
\end{equation}
where $\sigma_n$ is the volume of the $n$-dimensional unit ball.
\end{corollary}
\begin{proof}
For a ball $B$ of radius $r$ in $\R^n$ centered at the origin we have
$$
\frac{\int_B \vert x\vert^2\:dx}{\vol(B)}=\frac{nr^2}{n+2}
$$
If $\Delta$ is a ball of radius $r$ centered at the origin,
then $\mathfrak N(\pi)$ is a ball of radius $r$ in $\mathfrak k^*$.
Therefore, from (\ref{eqAsymp1}) it follows that
$r^2_\Delta=(n+2)^{-1}(\alpha,\alpha+2\rho)^{-1}r^2$.
Now the desired equality follows from Theorem \ref{thmAsymp},
since the volume of a ball of radius $r_\Delta$ equals
the right-hand side of equality (\ref{volAsymp}).
\end{proof}
\emph{Proof of Theorem \ref{thmAsymp1}.}
According to (\ref{volAsymp}),
the desired statement follows from Theorem \ref{thmMixed}.
\section{Probability of the root being real}\label{proportion}
\subsection{BKK Theorem for Reductive Groups}\label{BKK}
Let $K^\C$ be the complexification of a compact group $K$.
It is a complex connected $n$-dimensional reductive Lie group
such that $K$ is a maximal compact subgroup of $K^\C$.
We consider finite-dimensional holomorphic representations $\mu_1,\ldots,\mu_n$ of $K^\C$,
and the complex vector spaces of $\mu_i$-polynomials ${\rm Trig}(\mu_i)$
(recall that a $\mu_i$-polynomial is a linear combination of matrix elements of the representation $\mu_i$).
To any system of $n$ nonzero $\mu_i$-polynomials $f_i\in{\rm Trig}(\mu_i)$
we associate the point
$\iota(f_1,\ldots,f_n)=(\C f_1)\times\ldots\times(\C f_n)\in\P_{1,\C}\times\ldots\times\P_{n,\C}$,
where $\P_{i,\C}$ is the complex projective space whose points are
one-dimensional subspaces of ${\rm Trig}(\mu_i)$.
We use the following standard result from algebraic geometry.
\begin{proposition}\label{prConstAl}
There exist a number $N(\mu_1,\ldots,\mu_n)$ and an algebraic hypersurface $H$ in $\P_{1,\C}\times\ldots\times\P_{n,\C}$
such that the following holds.
For any $n$ $\mu_i$-polynomials $f_i\in{\rm Trig}(\mu_i)$ with $\iota(f_1,\ldots,f_n)\not\in H$,
the number of their common zeros equals $N(\mu_1,\ldots,\mu_n)$.
\end{proposition}
Below we provide a known geometric formula for $N(\mu_1,\ldots,\mu_n)$ (Theorem \ref{thmRed1})
and a variant of this formula (Theorem \ref{thmRed}),
which is used to derive Theorem \ref{thmAntiKac2}.

Consider the decomposition
$$
\mu=\bigoplus_{\lambda\in\Lambda\subset\Z^k\cap\mathfrak C^*,\:0<m_\lambda\in\Z} m_\lambda\:\mu_\lambda
$$
of a representation $\mu$ into a sum of irreducible representations $\mu_\lambda$
with highest weights $\lambda$ and multiplicities $m_\lambda$.
\begin{definition}\label{dfweightedpolyhedron}
We denote the Weyl group orbit of a point $\lambda\in\mathfrak t^*$ by $W^*(\lambda)$.
The compact convex set
$$
\Delta(\mu)={\rm conv}\left(\bigcup_{\lambda\in\Lambda}W^*(\lambda)\right)
$$
is called the \emph{weight polytope of the representation} $\mu$.
Let $\mathfrak N(\mu)\subset \mathfrak k^*$ denote the union of coadjoint orbits of $K$
intersecting the weight polytope $\Delta(\mu)$
(we identify $\mathfrak t^*$ with the set of fixed points of the coadjoint action of $T^k$ on $\mathfrak k^*$).
In what follows, we call $\mathfrak N(\mu)$ the \emph{Newton body of the representation} $\mu$.
\end{definition}
\begin{corollary}\label{corNewtTor}
{\rm(1)} Let $\mu_T$ be a restriction of the representation $\mu$ to the maximal torus $T^k$ in $K$,
and let $\Lambda_T\subset \mathfrak t^*$ be the set of weights of the torus representation $\mu_T$.
Then $\Delta(\mu)={\rm conv}(\Lambda_T)$.

{\rm (2)}
The set $\mathfrak N(\mu)$ is convex.

{\rm (3)}
$\Delta(\mu)=\pi(\mathfrak N(\mu))$,
where $\pi$ is the projection $\mathfrak k^*\to \mathfrak t^*$.

{\rm (4)}
$\mathfrak N(\mu\otimes \pi)=\mathfrak N(\mu)+ \mathfrak N(\pi)$.
\end{corollary}
\begin{proof}
Statement (1) follows from the theory of highest weights.
It is known that for any $\zeta\in\mathfrak t^*$,
the projection $\pi$ of the coadjoint orbit ${\rm Ad}(K)(\zeta)$ onto $\mathfrak t^*$
coincides with the convex hull of the Weyl group orbit $W^*\zeta$; see \cite{K}.
From this, statements (2) and (3) follow.
Statement (4) follows from standard properties of representation weights.
\end{proof}
Let $F$ be a homogeneous polynomial of degree $p$ on $\mathfrak t^*$.
For a convex polytope $\Delta\subset\mathfrak t^*$, set
\begin{equation}\label{eqv}
  I(\Delta;F)=\int_\Delta F\: d\upsilon,
\end{equation}
where the $W^*$-invariant measure $\upsilon$ on $\mathfrak t^*$ is chosen
so that the volume of the fundamental parallelepiped of the character lattice $\Z^k$ equals $1$.
It is known (see, e.g., \cite{AKK}) that
$I(\Delta;F)$ is a homogeneous polynomial of degree $k+p$
on the space of virtual convex polytopes in $\mathfrak t^*$.
Denote by $J(\Delta_1,\ldots,\Delta_{k+p};F)$
its polarization, i.e., the symmetric $(k+p)$-linear form on the space of virtual convex polytopes
such that $J(\Delta_1,\ldots,\Delta_{k+p};F)=I(\Delta;F)$
when $\Delta_1=\ldots=\Delta_{k+p}=\Delta$.

First, we state the reductive version of the BKK theorem from \cite{AKK}.
According to the \emph{Weyl dimension formula},
the dimension of a representation $\mu_\lambda$ with highest weight $\lambda$ is $F_W(\lambda)$, where
$F_W$ is a certain polynomial on $\mathfrak t^*$ of degree $(2n-k)/2$.
We denote the highest homogeneous component of $F_W$ by $\phi$.
\begin{theorem}\label{thmRed1}
For finite-dimensional holomorphic representations $\mu_1,\ldots,\mu_n$ of $K^\C$, we have
$$
N(\mu_1,\ldots,\mu_n)=\frac{n!}{\vert W\vert} J(\Delta(\mu_1),\ldots,\Delta(\mu_n);\phi^2),
$$
where $\Delta(\mu_i)$ is the weight polytope of the representation $\mu_i$ (see Definition \ref{dfweightedpolyhedron}).
\end{theorem}
Next, we use the Weyl integration formula;
see (\ref{WeyIntFormula}) in \textsection\ref{sAsymp}.
Recall its statement:
for a $K$-coadjoint-invariant function $f\colon\mathfrak k^*\to\R$,
$$
 \int_{\mathfrak k^*}f \:d\nu=\frac{\vol_\tau(K)}{\vert W\vert\vol_\tau(T^k)}\int_{\mathfrak t^*}P^2(\lambda)\:f\:d\nu,
$$
where $P(\lambda)=\prod_{\theta\in R^+}(\lambda,\theta)$.
\begin{theorem}\label{thmRed}
We have
$$
N(\mu_1,\ldots,\mu_n)=\frac{n!}{P^2(\rho)}\:\vol_\tau\left(\mathfrak N(\mu_1),\ldots,\mathfrak N(\mu_n)\right),
$$
where $\rho=\frac{1}{2}\sum_{\beta\in R^+}\beta$.
\end{theorem}
\begin{proof}
Consider the invariant metric $\tau^s=s\tau$ on $\mathfrak k$,
where $s^{-k}$ equals the volume of the fundamental parallelepiped of the character lattice in $\mathfrak t^*$
measured using $\tau$ (see the beginning of \textsection\ref{s2.1}).
Let $\upsilon$ be the corresponding Lebesgue measure on $\mathfrak t^*$.
Then, for the measure $\upsilon$ and the functional $I$ from (\ref{eqv}), the condition of Theorem \ref{thmRed1} is satisfied.

Recall that
according to the Weyl formula for the dimension of an irreducible representation with highest weight $\lambda$,
$\dim(\mu_\lambda)=P(\lambda+\rho)/P(\rho)$.
Hence
$$F_W(\lambda)=\frac{P(\lambda+\rho)}{P(\rho)},\quad \phi^2=\frac{P^2(\lambda+\rho)}{P^2(\rho)},\quad I(\Delta;\phi^2)=\int_\Delta\frac{P^2(\lambda+\rho)}{P^2(\rho)} \: d\upsilon.$$
Now, Theorem \ref{thmRed1} implies
$$
N(\mu,\ldots,\mu)=\frac{n!}{\vert W\vert\:P^2(\rho)} \int_{\Delta(\mu)}P^2(\lambda)d\upsilon.
$$
Let $\nu^s$ denote the Lebesgue measure in $\mathfrak k^*$ corresponding to $\tau^s$.
By definition, $\nu^s=s^n\nu$.
Applying the Weyl integration formula to $\nu^s$ and to the function $f=\vartheta P^2(\lambda)$,
where $\vartheta$ is the characteristic function of the weight polytope $\Delta(\mu)$,
and noting that $\vol_{\tau^s}(T^k)=1$ and $\vol_{\tau^s}(K)=s^n$,
we obtain
\begin{multline*}
N(\mu,\ldots,\mu)=\frac{n!}{\vert W\vert P^2(\rho)} \int_{\mathfrak t^*}\vartheta\:P^2(\lambda)\:d\nu^s
=\\
\frac{n!}{s^n\:P^2(\rho)}\int_{\Delta(\mu)}P^2(\lambda)\:d\nu^s
=
\frac{n!}{P^2(\rho)}\vol_\tau(\mathfrak N(\mu)).
\end{multline*}
Applying the polarization formula to the homogeneous polynomial $I(\Delta;\frac{P^2(\lambda+\rho)}{P^2(\rho)})$ in the argument $\Delta$,
and also to the volume polynomial in the space of virtual compact sets in $\mathfrak k^*$ (see, e.g., \cite{KK}),
we obtain the desired statement.
\end{proof}
\subsection{Probability of the root being real}\label{expRoot}
\begin{theorem}\label{thmProportion}
Let $\pi_1,\ldots,\pi_n$ be real representations of a simple Lie group $K$,
and let $\pi^\C_i$ be the complexification of $\pi_i$.
Denote by $\mathcal P(\pi_1,\ldots,\pi_n)$ the probability that a root of a system of $n$
random $\pi_i$-polynomials $f_1=\ldots=f_n=0$
is real.
Then
\begin{equation}\label{eqCR}
\mathcal P(\pi_1,\ldots,\pi_n)=\frac{P^2(\rho)}{(2\pi)^n}\:\frac{\vol\left({\rm Ell}(\pi_1),\ldots,{\rm Ell}(\pi_n)\right)}{\vol\left(\mathfrak N(\pi^\C_1),\ldots,\mathfrak N(\pi^\C_n)\right)},
\end{equation}
where $\rho=\frac{1}{2}\sum_{\beta\in R^+}\beta$ is the half-sum of positive roots of $K$,
and $P(\lambda)=\prod_{\beta\in R^+}(\lambda,\beta)$.
\end{theorem}
\begin{proof}
Then, for any $i\le n$, the real vector subspace ${\rm Trig}(\pi_i)$
is dense in the complex vector space ${\rm Trig}(\pi^\C_i)$ with respect to the Zariski topology.
Let $\P_{i}$ and $\P_{i,\C}$ be the real and complex projectivizations of ${\rm Trig}(\pi_i)$ and ${\rm Trig}(\pi^\C_i)$, respectively.
Consider the embedding
$$
\iota\colon\P_1\times\ldots\times\P_n\to\P_{1,\C}\times\ldots\times\P_{n,\C}.
$$
The image of $\iota$ is Zariski dense in $\P_{1,\C}\times\ldots\times\P_{n,\C}$.
Therefore, by Proposition \ref{prConstAl}, $\iota^{-1}H$ is contained in a closed real hypersurface in $\P_1\times\ldots\times\P_n$.
It follows that the number of roots of real $\pi_i$-polynomials almost always equals $N(\pi^\C_1,\ldots,\pi^\C_n)$.
Hence
$$
  \mathcal P(\pi_1,\ldots,\pi_n)=\frac{\mathfrak M(\pi_1,\ldots,\pi_n)}{N(\pi^\C_1,\ldots,\pi^\C_n)}.
$$
The desired statement follows from Theorems \ref{thmMixed} and \ref{thmRed}.
\end{proof}
\subsection{Asymptotic probability of the root being real}\label{AsympExpReal}
Let $B$ be a compact, convex, centrally symmetric subset of $\mathfrak t^*$.
Assume $B$ is invariant under the Weyl group $W^*$.
Then the set $\Lambda(B)=B\cap\mathfrak C^*\cap\Z^k$ is symmetric; see Definition \ref{dfSymm}.
Consider a sequence of sets $\Lambda(mB)$ and the corresponding sequence of representations
\begin{equation}\label{eqpim}
  \pi_m=\bigoplus_{(\lambda,\lambda')\in\Lambda'(mB)}\pi_{\lambda,\lambda'}.
\end{equation}
We further assume $B$ is a unit ball
and consider the holomorphic extensions of $\pi_m$-polynomials to polynomials in ${\rm Trig}(K^\C)$ as analogues of real Laurent polynomials of degree $m$.
Denote by ${\mathcal P}(\pi)$
the probability that a root of a system of $n$ $\pi$-polynomials is real,
i.e., ${\mathcal P}(\pi)={\mathcal P}(\pi,\ldots,\pi)$.
\begin{theorem}\label{thmAntiKac2}
Let $K$ be a simple group.
Then
$$
\lim_{m\to\infty} {\mathcal P}(\pi_m)=
\frac{P^2(\rho)}{(2\pi)^n(n+2)^{n/2}(\alpha,\alpha+2\rho)^{n/2}},
$$
where $\alpha$ is the highest root of $K$ (i.e., the highest weight of the adjoint representation $\mu_\alpha$).
\end{theorem}
\begin{proof}
From (\ref{eqCR}) we have
$$
\mathcal P(\pi_m)=\frac{P^2(\rho)}{(2\pi)^n}\:\frac{\vol\left({\rm Ell}(\pi_m)\right)}{\vol\left(\mathfrak N(\pi_m^\C)\right)}.
$$
By Definition \ref{dfweightedpolyhedron},
the Newton body $\mathfrak N(\pi^\C_m)$ is the union of coadjoint $K$-orbits
through points in $\rm{conv}((mB)\cap\Z^k)$.
Hence $\mathfrak N(\pi^\C_m)$ asymptotically coincides with a ball of radius $m$ in $\mathfrak k^*$.
Moreover,
if $K$ is simple, the Newton ellipsoid is also a ball; see Corollary \ref{Ball}.
Thus, computing the limiting probability $\lim_{m\to\infty} {\mathcal P}(\pi_m)$
reduces to
computing the asymptotic radius of the ellipsoid ${\rm Ell}(\pi_m)$ as $m\to\infty$.
Applying estimate (\ref{volAsymp}) for $\vol({\rm Ell}(\pi_m))$ from Corollary \ref{corAsymp1},
we obtain the desired statement.
\end{proof}
\begin{remark}
Using \cite{Mcd,Hash,HCh},
the equality in Theorem \ref{thmAntiKac2} can be expressed in a more topological form.
However, this topological expression appears more cumbersome.
\end{remark}

\end{document}